\documentclass[11pt]{amsart}
\usepackage{amsmath, amsthm, amscd, amsfonts, amssymb, graphicx, color,float,pgf,tikz}
\usepackage{amssymb,fontenc}
\usepackage{latexsym,wasysym,mathrsfs}
\usepackage{hyperref}
\usetikzlibrary{arrows}
\makeatletter \oddsidemargin.9375in \evensidemargin \oddsidemargin
\newtheorem{theorem}{Theorem}[section]

\theoremstyle{definition}

\theoremstyle{remark}
\newtheorem{remark}[theorem]{Remark}
\numberwithin{equation}{section}

\newcommand{\be}{\begin{equation}}
\newcommand{\ee}{\end{equation}}

\numberwithin{equation}{section}

\begin{document}

\noindent {\tiny Sahand Communications in Mathematical Analysis (SCMA) Vol. $\cdots$ No. $\cdots$(2014), $\cdots$-$\cdots$\\
\url{http://scma.maragheh.ac.ir}}   \\[0.50in]


\title[character space of BSE-functions]{On character space of the algebra of BSE-functions}

\author[M. Fozouni]{Mohammad Fozouni$^1$$^{*}$}

\address{ $^{1}$ Department of Mathematics and Statistics, Faculty of Basic Sciences and Engineering,  Gonbad Kavous University, P.O.Box 163, Gonbad Kavous, Iran.}

\email{fozouni@hotmail.com}




\keywords{Banach algebra, BSE-function, character space, locally compact group.\\
\indent Received: dd mmmm yyyy,    Accepted: dd mmmm yyyy.
\\
\indent $^{*}$ Corresponding author}
\maketitle
\hrule width \hsize \kern 1mm


\begin{abstract}
Suppose that $A$ is a semi-simple and commutative Banach algebra. In this paper we try to characterize the character space of the Banach algebra $C_{\rm{BSE}}(\Delta(A))$ consisting of all  BSE-functions on $\Delta(A)$ where $\Delta(A)$ denotes the character space of $A$. Indeed, in the case that $A=C_0(X)$ where $X$ is a non-empty locally compact Hausdroff space, we give a complete characterization of $\Delta(C_{\rm{BSE}}(\Delta(A)))$ and in the general case we give a partial answer.

 Also, using the Fourier algebra, we show that $C_{\rm{BSE}}(\Delta(A))$ is not a $C^*$-algebra in general.
Finally for some subsets $E$ of $A^*$, we define the subspace of BSE-like functions on $\Delta(A)\cup E$ and give a nice application of this space related to Goldstine's theorem.
\end{abstract}
\maketitle
\vspace{0.1in}
\hrule width \hsize \kern 1mm

\section{Introduction and Preliminaries}
Suppose that $A$ is a semi-simple commutative Banach algebra  and $\Delta(A)$ is the character space of $A$, i.e., the space of all non-zero homomorphisms from $A$ into $\mathbb{C}$.

 A bounded continuous function $\sigma$ on $\Delta(A)$ is called a BSE-function if there exists a constant $C>0$ such that for each $\phi_{1},...,\phi_{n}\in\Delta(A)$ and complex numbers $c_{1},...,c_{n}$, the inequality
\begin{equation*}
\left|\sum_{i=1}^{n}c_{i}\sigma(\phi_{i})\right|\leq C\left\|\sum_{i=1}^{n}c_{i}\phi_{i}\right\|_{A^{*}}
\end{equation*}
holds. For each $\sigma\in C_{\mathrm{BSE}}(\Delta(A))$ we denote by $\|\sigma\|_{\mathrm{BSE}}$ the infimum of such $C$. Let $C_{\mathrm{BSE}}(\Delta(A))$ be the set of all $\mathrm{BSE}$-functions. We have a good characterization of $C_{\mathrm{BSE}}(\Delta(A))$ as follows:
\begin{theorem}\label{1}
$C_{\mathrm{BSE}}(\Delta(A))$ is equal to the set of all $\sigma\in C_{b}(\Delta(A))$ for which there exists a bounded net $\{x_{\lambda}\}$ in $A$ with $\lim_{\lambda}\phi(x_{\lambda})=\sigma(\phi)$ for all $\phi\in \Delta(A)$
\end{theorem}
\begin{proof}
 see \cite[Theorem 4 (i)]{Takahasi}.
\end{proof}

Using the above characterization one can see that $C_{\rm{BSE}}(\Delta(A))$ is unital if and only if $A$ has a bounded weak approximate identity  in the sense of Lahr and Jones. We recall that a net $\{x_{\alpha}\}$ in $A$ is called a bounded weak approximate identity (b.w.a.i) for $A$ if $\{x_{\alpha}\}$ is bounded in $A$ and
\begin{equation*}
\lim_{\alpha}\phi(x_{\alpha}a)=\phi(a)\hspace{0.5cm}(\phi\in \Delta(A), a\in A),
\end{equation*}
or equivalently, $\lim_{\alpha}\phi(x_{\alpha})=1$ for each $\phi\in \Delta(A)$.

Also, Theorem \ref{1}, gives the following definition of $\|\cdot\|_{\rm{BSE}}$:
\begin{align*}
  \|\sigma\|_{\rm{BSE}}=\inf\{ & \beta>0 : \exists \{x_{\lambda}\} \text{ in A with } \|x_{\lambda}\|\leq \beta, \\
   & \lim_{\lambda}\phi(x_{\lambda})=\sigma(\phi)\ (\phi\in \Delta(A)) \}.
\end{align*}

The theory of  BSE-algebras for the first time introduced and investigated by Takahasi and Hatori; see \cite{Takahasi} and two other notable works \cite{Kaniuth, Kamali2}. In \cite{Kamali2}, the authors answered to a question raised in \cite{Takahasi}. Examples of BSE-algebras are the group algebra $L^1(G)$ of a locally compact abelian
group $G$, the Fourier algebra $A(G)$ of a locally compact amenable
group $G$, all commutative $C^*$-algebras, the disk algebra, and the Hardy algebra on the
open unit disk. We recall that a commutative and without order Banach algebra $A$ is a type I-BSE algebra if $\widehat{M(A)}=C_{\rm{BSE}}(\Delta(A))=C_b(\Delta(A))$, where $M(A)$ denotes the multiplier algebra of $A$ and $\widehat{M(A)}$ denotes the space of all $\widehat{T}$ which defined by $\widehat{T}(\varphi)\varphi(x)=\widehat{T(x)}(\varphi)$ for all $\varphi\in \Delta(A)$. Note that $x\in A$ should satisfies $\varphi(x)\neq 0$.

In this paper,  we give a partial characterization of the character space of $C_{\rm{BSE}}(\Delta(A))$ where $A$ is a semi-simple commutative Banach algebra. Indeed, we show that if  $A$ has a b.w.a.i and $C_{\rm{BSE}}(\Delta(A))$ is an ideal in $C_b(\Delta(A))$, then
$$\Delta(C_{\rm{BSE}}(\Delta(A)))=\overline{\Delta(A)}^{w^*}.$$

 Also, we give a negative answer to this question;\\
    Whether $(C_{\rm{BSE}}(\Delta(A)), \|\cdot\|_{\rm{BSE}})$ is a $C^*$-algebra? At the final section of this paper we study the space of BSE-like functions on subsets of $A^*$ which containing $\Delta(A)$ and as an application of this space we give a nice relation with Goldstine's theorem.

\section{Character space of $C_{\mathrm{BSE}}(\Delta(A))$}

In view of \cite[Lemma 1]{Takahasi}, $C_{\rm{BSE}}(\Delta(A))$ is a semi-simple commutative Banach algebra. So, the character space of $C_{\rm{BSE}}(\Delta(A))$ should be non-empty and one may ask: Is there a characterization of $\Delta(C_{\rm{BSE}}(\Delta(A)))$ for an arbitrary Banach algebra $A$?

 In the sequel of this section we give a partial answer to this question.
Let  $X$ be a non-empty locally compact Hausdroff space and put
$$C_{\mathrm{BSE}}(X):=C_{\mathrm{BSE}}(\Delta(C_{0}(X))).$$

To proceed further we recall some notions. Let $X$ be a non-empty locally compact Hausdroff space. A function algebra (FA) on $X$ is a subalgebra $A$ of $C_b(X)$ that seperates strongly the points of $X$, that is, for each $x, y\in X$ with $x\neq y$, there exists $f\in A$ with $f(x)\neq f(y)$ and for each $x\in X$, there exists $f\in A$ with $f(x)\neq 0$.
A Banach function algebra (BFA) on $X$ is a function algebra $A$ on $X$ with a norm $\|\cdot\|$ such that $(A, \|\cdot\|)$ is a Banach algebra. 

A topological space $X$ is completely regular if every non-empty closed set and every singleton disjoint from it can be separated by continuous functions.

\begin{theorem} Let  $X$ be a non-empty locally compact Hausdroff space. Then
  $C_{\mathrm{BSE}}(X)$ is a unital \rm{BFA} and its character space is homeomorphic to $\overline{X}^{w^{*}}$, that is, $\Delta(C_{\mathrm{BSE}}(X))=\overline{\{\phi_{x} : x\in X\}}^{w^*}$.
\end{theorem}
\begin{proof}
By \cite[Lemma 1]{Takahasi}, $C_{\mathrm{BSE}}(X)$ is a subalgebra of $C_{b}(X)$ and $\|\cdot\|_{\rm{BSE}}$ is a complete algebra norm.
Since $C_{0}(X)$ has a bounded approximate identity,  $C_{\rm{BSE}}(X)$ is unital. So, for each $x\in X$, there exists $f\in C_{\mathrm{BSE}}(X)$ with $f(x)\neq 0$. On the other hand, using the Urysohn lemma, for each $x, y\in X$ with $x\neq y$ one can see that there exists $f\in C_{\mathrm{BSE}}(X)$ such that $f(x)\neq f(y)$.

Finally, since $X$ is a locally compact Hausdroff space, it is  completely regular by \cite[Corollary 2.74]{AK}.
  On the other hand, by \cite[Theorem 3]{Takahasi}, we know that $C_{0}(X)$ is a type I-BSE algebra. Therefore,  $C_{\textrm{BSE}}(X)=C_{b}(\Delta(C_{0}(X)))=C_{b}(X)$.
Also,  for every $f\in C_{\rm{BSE}}(X)$, by the remark after Theorem 4 of \cite{Takahasi}, we have $\|f\|_{X}\leq \|f\|_{\rm{BSE}}$.
Also, by the Open mapping theorem there exists a positive constant $M$ such that $\|f\|_{\rm{BSE}}\leq M\|f\|_X$. So, $C_{\rm{BSE}}(X)$ and $C_b(X)$ are topologically isomorphic, and so $\Delta(C_{\textrm{BSE}}(X))$ and $\Delta(C_{b}(X))$ are homeomorphic.
Now, by using \cite[Theorem 2.4.12]{Kaniuth2}, we have
 $$\Delta(C_{\textrm{BSE}}(X))=\Delta(C_{b}(X))=\overline{X}^{w^{*}}=\overline{\{\phi_{x} : x\in X\}}^{w^{*}}.$$
\end{proof}

\begin{remark}\label{Rem1}
In general for a commutative Banach algebra $A$, we have the following conditions concerning the character space of $C_{\rm{BSE}}(\Delta(A))$:
\begin{enumerate}
  \item If $C_{\rm{BSE}}(\Delta(A))=C_{b}(\Delta(A))$, then
  $$\Delta(C_{\rm{BSE}}(\Delta(A)))=\overline{\Delta(A)}^{w^{*}}.$$
  Examples of Banach algebras $A$ satisfying $C_{\rm{BSE}}(\Delta(A))=C_{b}(\Delta(A))$ are  finite dimensional Banach algebras and commutative $C^*$-algebras; see the remark on page 609 of \cite{TH2}. Also, see \cite[Lemma 2]{Takahasi} for a characterization of Banach algebras $A$ for which satisfying  $C_{\rm{BSE}}(\Delta(A))=C_{b}(\Delta(A))$.
  \item If  $A$ has a b.w.a.i, then $C_{\rm{BSE}}(\Delta(A))$ is unital and so $\Delta(C_{\rm{BSE}}(\Delta(A)))$ is compact and hence it is $w^*$-closed.
     On the other hand, we know that $$\Delta(A)\subseteq \Delta(C_{\rm{BSE}}(\Delta(A))),$$
     in the sense that for each $\varphi\in \Delta(A)$, $f_{\varphi}: C_{\rm{BSE}}(\Delta(A))\longrightarrow \mathbb{C}$ defined by $f_{\varphi}(\sigma)=\sigma(\varphi)$ is an element of $\Delta(C_{\rm{BSE}}(\Delta(A)))$. Note that $f_\varphi\neq 0$, since in this case $C_{\rm{BSE}}(\Delta(A))$ is unital and $f_\varphi(1)=1$.
     So $$\overline{\Delta(A)}^{w^{*}}\subseteq \overline{\Delta(C_{\rm{BSE}}(\Delta(A)))}^{w^*}=\Delta(C_{\rm{BSE}}(\Delta(A))).$$
    \item On can see that if $(B, \|\cdot\|_B)$ is a Banach algebra which containing the Banach algebra $(C,\|\cdot\|_C)$ as a two-sided ideal, then every $\varphi\in \Delta(C)$ extends to one $\widetilde{\varphi}\in \Delta(B)$. Now, let $C=C_{\rm{BSE}}(\Delta(A))$ and $B=C_{b}(\Delta(A))$. If $C$ is an ideal in $B$, then
         $$\Delta(C_{\rm{BSE}}(\Delta(A)))=\Delta(C)\subseteq \Delta(B)=\overline{\Delta(A)}^{w^*}.$$
    \item Suppose that $B$ is a commutative semi-simple Banach algebra such that $\Delta(B)$ is compact. Then $B$ is unital; see \cite[Theorem 3.5.5]{Kaniuth2}.
Now, If $A$ has no b.w.a.i, then $$\Delta(C_{\rm{BSE}}(\Delta(A)))\neq \overline{\Delta(A)}^{w^{*}}.$$
     Because if $\Delta(C_{\rm{BSE}}(\Delta(A)))=\overline{\Delta(A)}^{w^{*}},$  by using the above assertion, $C_{\rm{BSE}}(\Delta(A))$ is unital,  since $\overline{\Delta(A)}^{w^{*}}=\Delta(C_{b}(\Delta(A)))$ is compact and  $C_{\rm{BSE}}(\Delta(A))$ is a semi-simple commutative Banach algebra. Therefore, $A$ has a b.w.a.i which is impossible.

     \item If $A$ has a b.w.a.i and $C_{\rm{BSE}}(\Delta(A))$ is an ideal of $C_b(\Delta(A))$, then using parts (ii) and (iii), we have
         $$\Delta(C_{\rm{BSE}}(\Delta(A)))=\overline{\Delta(A)}^{w^*}.$$
\end{enumerate}
\end{remark}
\section{$\left(C_{\rm{BSE}}(\Delta(A)), \|\cdot \|_{\rm{BSE}}\right)$ is not a $C^{*}$-algebra}

The theory of $C^*$-algebras is very fruitful and applied. As an advantage of this theory, especially in Harmonic Analysis, one can see the $C^*$-algebra approach for defining a locally compact quantum group; see \cite{KV00}. So,  verifying the Banach algebras from a $C^*$-algebraic point of view is very helpful. In this section, using a result due to Kaniuth and \"{U}lger in \cite{Kaniuth}, we show that $\left(C_{\rm{BSE}}(\Delta(A)), \|\cdot \|_{\rm{BSE}}\right)$ is not a $C^{*}$-algebra in general. On the other hand, there is a question which left open that, under what conditions on $A$, $\left(C_{\rm{BSE}}(\Delta(A)), \|\cdot \|_{\rm{BSE}}\right)$ is  a $C^{*}$-algebra?

In the sequel for each locally compact group $G$, let $A(G)$ denote the Fourier algebra and $B(G)$ denote the Fourier-Stieltjes algebra introduced by Eymard; see \cite[\S 19]{pier}. Also, let $\widehat{G}$ denote the dual group of $G$ and $M(G)$ denote the Measure algebra; see \cite[\S 3.3]{Dales}. For the convenience of reader we give the definitions of $A(G)$ and $B(G)$ as follows:

Let $G$ be a locally compact group.  Suppose that $A(G)$ denotes the subspace of $C_{0}(G)$ consisting of functions of the form $u=\sum_{i=1}^{\infty}f_{i}\ast \widetilde{g_{i}}$ where $f_{i}, g_i\in L^{2}(G)$, $\sum_{i=1}^{\infty}||f_{i}||_{2}||g_{i}||_{2}<\infty$ and $\widetilde{f}(x)=\overline{f(x^{-1})}$ for all $x\in G$. The space $A(G)$ with the pointwise operation and the following norm is a Banach algebra,
\begin{equation*}
||u||_{A(G)}=\inf\{\sum_{i=1}^{\infty}||f_{i}||_{2}||g_{i}||_{2}: u=\sum_{i=1}^{\infty}f_{i}\ast \widetilde{g_{i}}\},
\end{equation*}
which we call it the Fourier algebra. It is obvious that for each $u\in A(G)$, $||u||\leq ||u||_{A(G)}$ where $||u||$ is the norm of $u$ in $C_{0}(G)$.

Now let $\Sigma$ denote the equivalence class of all irreducible representations of $G$. Then $B(G)$ consisting of all functions $\phi$ of the form $\phi(x)=<\pi(x)\xi, \eta>$ where $\pi\in \Sigma$ and $\xi, \eta$ are elements of $H_{\pi}$, the Hilbert space associated to the representation $\pi$. It is well-known that $A(G)$ is a closed ideal of $B(G)$.

Also, recall that an involutive Banach algebra $A$ is called a $C^*$-algebra if its norm satisfies $\|aa^*\|=\|a\|^2$ for each $a\in A$. We refer the reader to \cite{Murphy} to see a complete description of $C^*$-algebras.

In the following remark we give the main result of this section.

\begin{remark} In general $\left(C_{\rm{BSE}}(\Delta(A)), \|\cdot \|_{\rm{BSE}}\right)$ is not a $C^{*}$-algebra, that is, there is not any involution "*" on $C_{\rm{BSE}}(\Delta(A))$ such that
  $$\|\sigma^*\sigma \|_{\rm{BSE}}=\|\sigma \|_{\rm{BSE}}^2 \ \ \ \forall \sigma\in C_{\rm{BSE}}(\Delta(A)).$$

  Because we know that every commutative $C^{*}$-algebra is a BSE-algebra. For a non-compact locally compact Abelian group $G$ take $A=A(G)$. By \cite[Theorem 5.1]{Kaniuth}, we know that $C_{\rm{BSE}}(\Delta(A))=B(G)$ and for each $u\in B(G)$, $\|u\|_{B(G)}=\|u\|_{\rm{BSE}}$. But $B(G)=M(\widehat{G})$ and it is shown in \cite{Takahasi} that $M(\widehat{G})$  and hence $B(G)$ is not a BSE-algebra. Therefore, $\left(C_{\rm{BSE}}(\Delta(A)), \|\cdot \|_{\rm{BSE}}\right)$ is not a $C^{*}$-algebra.

  As the second example, let $G$ be a locally compact Abelian  group. It is well-known that $C_{\rm{BSE}}(\Delta(L^1(G)))$ is isometrically isomorphic to $M(G)$, where $L^1(G)$ denotes the group algebra; see the last remark on page 151 of \cite{Takahasi}. On the other hand, by the Gelfand-Nimark theorem we know that every commutative $C^*$-algebra should be symmetric. But in general $M(G)$ is not symmetric, i.e., the formula $\widehat{\mu^*}(\xi)=\overline{\widehat{\mu}}(\xi)$ does not hold for every $\xi\in \Delta(M(G))$. For example if $G$ is non-discrete then by \cite[Theorem 5.3.4]{Rudin}, $M(G)$ is not symmetric and hence fails to be a $C^*$-algebra.
\end{remark}
It is a good question to characterize Banach algebras $A$ for which $\left(C_{\rm{BSE}}(\Delta(A)), \|\cdot \|_{\rm{BSE}}\right)$ is a $C^*$-algebra.
\section{BSE-functions on subsets of $A^*$}

Suppose that $A$ is a Banach algebra and $E\subseteq A^*\setminus \Delta(A)$. A complex-valued bounded continuous function $\sigma$ on $\Delta(A)\cup E$ is called a BSE-like function if there exists an $M>0$ such that for each $f_1, f_2, f_3, \ldots f_n\in \Delta(A)\cup E$ and complex numbers $c_1, c_2, c_3, \ldots, c_n$, 
\begin{equation}\label{eq1}
\left|\sum_{i=1}^{n}c_{i}\sigma(f_{i})\right|\leq M\left\|\sum_{i=1}^{n}c_{i}f_{i}\right\|_{A^*}.
\end{equation}
We show the set of all the BSE-like functions on $\Delta(A)\cup E$ by $C_{\rm{BSE}}(\Delta(A), E)$ and let $\|\sigma\|_{\rm{BSE}}$ be the infimum of all $M$ satisfying relation \ref{eq1}. Obviously, $C_{\rm{BSE}}(\Delta(A), E)$ is a linear subspace of $C_b(\Delta(A)\cup E)$ and we have
 $$\{\sigma_{|\Delta(A)} : \sigma\in C_{\rm{BSE}}(\Delta(A), E)\}\subseteq C_{\rm{BSE}}(\Delta(A)).$$

Clearly, $\iota_{A}(A)\subseteq C_{\rm{BSE}}(\Delta(A), E)$ where $\iota_A: A\longrightarrow A^{**}$ is the natural embedding. For $a\in A$, we let $\widehat{a}=\iota_A(a)$ and $\widehat{A}=\iota_A(A)$.

To proceed further, we recall the Helly theorem.
\begin{theorem}\label{Helly}
  (\rm{Helly}) Let $(X, \|\cdot\|)$ be a normed linear space over $\mathbb{C}$ and let $M>0$.  Suppose that $x_{1}^{*},\ldots, x_{n}^{*}$ are in $X^{*}$ and $c_{1}^{*},\ldots, c_{n}^{*}$ are in $\mathbb{C}$. Then the following are equivalent:
  \begin{enumerate}
    \item for all $\epsilon>0$, there exists $x_{\epsilon}\in X$ such that $\|x_{\epsilon}\|\leq M+\epsilon$ and $x_{k}^{*}(x_{\epsilon})=c_{k}$ for $k=1,\ldots, n$.
    \item for all $a_{1},\ldots, a_{n}\in \mathbb{C}$, $$\left|\sum_{i=1}^{n}a_{i}c_{i}\right|\leq M\left\|\sum_{i=1}^{n}a_{i}x_{i}^{*}\right\|_{X^*}.$$
  \end{enumerate}
 \end{theorem}
\begin{proof}
 See \cite[Theorem 4.10.1]{Larsen2}.
\end{proof}
As an application of Helly's theorem, we give the following characterization which is  similar to \cite[Theorem 4 (i)]{Takahasi}.

\begin{theorem}\label{Th:2}
 $C_{\rm{BSE}}(\Delta(A), E)$ is equal to the set of all $\sigma\in C_{b}(\Delta(A)\cup E)$ for which there exists a bounded net $\{x_{\alpha}\}$ in $A$ with $\lim_{\alpha}f(x_{\alpha})=\sigma(f)$ for all $f\in \Delta(A)\cup E$.
\end{theorem}
\begin{proof}
Suppose that $\sigma\in C_{b}(\Delta(A)\cup E)$ is such that there exists $\beta<\infty$ and a net $\{x_{\alpha}\}\subseteq X$ with $\|x_{\alpha}\|<\beta$ for all $\alpha$ and $\lim_{\alpha}f(x_{\alpha})=\sigma(f)$ for all $f\in \Delta(A)\cup E$. Let $f_{1}, \ldots, f_{n}$ be in $\Delta(A)\cup E$ and $c_1, \ldots, c_n$ be complex numbers. Then we have
\begin{align*}
  \left|\sum_{i=1}^{n}c_{i}\sigma(f_{i})\right|\leq & \left|\sum_{i=1}^{n}c_{i}f_{i}(x_{\alpha})\right|+ \left|\sum_{i=1}^{n}c_{i}(f_{i}(x_{\alpha})-\sigma(f_{i}))\right|\\
   \leq & \beta\left\|\sum_{i=1}^{n}c_{i}f_{i}\right\|+\sum_{i=1}^{n}|c_{i}||f_{i}(x_{\alpha})-\sigma(f_{i})|
\end{align*}
Taking the limit with respect to $\alpha$, we conclude that $\sigma\in C_{\rm{BSE}}(\Delta(A), E)$.

Conversely, let $\sigma\in C_{\rm{BSE}}(\Delta(A), E)$. Suppose that $\Lambda$ is the net consisting of all finite subsets of $\Delta(A)\cup E$. By Helly's theorem, for each $\epsilon>0$ and $\lambda\in \Lambda$, there exists $x_{(\lambda,\epsilon)}\in A$ with $\|x_{(\lambda,\epsilon)}\|\leq \|\sigma\|_{\rm{BSE}}+\epsilon$ and $f(x_{(\lambda,\epsilon)})=\sigma(f)$ for all $f\in \lambda$. Clearly, $\{(\lambda, \epsilon) : \lambda\in \Lambda, \epsilon>0\}$ is a directed set with $(\lambda_1,\epsilon_1)\preceq (\lambda_2,\epsilon_2)$ iff $\lambda_1\subseteq \lambda_2$ and $\epsilon_1\leq \epsilon_2$.
Therefore, we have
$$\lim_{(\lambda,\epsilon)}f(x_{(\lambda,\epsilon)})=\sigma(f)\ \ \ \ (f\in \Delta(A)\cup E).$$
\end{proof}
\begin{remark}
As an application of Theorem \ref{Th:2}, if $E=A^*\setminus \Delta(A)$, then one can see that $\overline{\widehat{A}}^{w^*}=A^{**}$, i.e., we conclude Goldstine's theorem. That is, $\widehat{A}$ with the $w^*$-topology of $A^{**}$ is dense in $A^{**}$.
\end{remark}
\begin{remark}
  We say that $A$ has a b.w.a.i respect to $E$ if there exists a bounded net $\{x_{\alpha}\}$ in $A$ with
  $$\lim_{\alpha}f(x_{\alpha})=1\ \ \  (f\in \Delta(A)\cup E).$$
Using Theorem \ref{Th:2},one can check that  $1\in C_{\rm{BSE}}(\Delta(A), E)$  if and only if $A$ has a b.w.a.i respect to $E$.
\end{remark}

We conclude this section with the following question.\\
{\bf Question}
Is $C_{\rm{BSE}}(\Delta(A), E)$ a commutative and semi-simple Banach algebra? If it is what is its character space?
  



\bigskip
\noindent
{\bf Acknowledgment.}
The author wish to thank the referee for his$\backslash$her suggestions. The author partially supported by a grant from Gonbad Kavous University.

\bibliographystyle{amsplain}

\vspace{0.1in}
\hrule width \hsize \kern 1mm
\end{document}